\documentclass{amsart}
\usepackage{amssymb}
\usepackage{epsfig}
\usepackage{psfrag}

\theoremstyle{plain}

\theoremstyle{definition}

\theoremstyle{definition}

\numberwithin{equation}{section}

\newcommand{\abs}[1]{\lvert #1 \rvert}

\newcommand{\R}{\ensuremath{\mathbb{R}}}

\newcommand{\Q}{\ensuremath{\mathbb{Q}}}

\begin{document}

\setlength{\textheight}{8.5in}
\setlength{\textwidth}{6in}
\setlength{\oddsidemargin}{0.2in}
\setlength{\topmargin}{-.5in}
\setlength{\evensidemargin}{0.2in}

\flushleft
\centerline{\large{{\bf Table of Contents for the Handbook of Knot
Theory}}}

\

\centerline{William W. Menasco and Morwen B. Thistlethwaite, Editors}

\

\begin{enumerate}

\item Colin Adams, {\it Hyperbolic knots}
\item Joan S. Birman and Tara Brendle {\it Braids and knots}
\item John Etnyre, {\it Legendrian and transversal knots}
\item Cameron Gordon {\it  Dehn surgery}
\item Jim Hoste {\it The enumeration and classification of knots and
links}
\item Louis Kauffman  {\it Diagrammatic methods for invariants of knots
and links}
\item Charles Livingston  {\it A survey of classical knot concordance}
\item Marty Scharlemann  {\it Thin position  }
\item Lee Rudolph  {\it Knot theory of complex plane curves}
\item DeWit Sumners  {\it The topology of DNA}
\item Jeff Weeks  {\it Computation of hyperbolic structures in knot
theory}
\end{enumerate}

\vfill\eject

\title[Hyperbolic Knots]
      {Hyperbolic Knots}

\author{Colin Adams}
\address{Dept. of Mathematics\\
         Williams College\\
         Williamstown, MA 01267}
\email{Colin.Adams@williams.edu}

\thanks{Thanks to Eric Schoenfeld for help with the figures.} 
\date{July, 2003}

\maketitle

\section{Introduction} \label{S:Introduction}

In 1978, Thurston revolutionized low dimensional topology when he demonstrated that many
 3-manifolds had hyperbolic metrics, or decomposed into pieces, many of which had 
hyperbolic metrics upon them. In view of the Mostow Rigidity 
theorem(~\cite{gM73}), when the volume
 associated with the manifold is finite, these hyperbolic metrics are unique. Hence
 geometric invariants coming out of the hyperbolic structure can be utilized to 
potentially distinguish between manifolds.

    One can either think of a hyperbolic 3-manifold as having a Riemannian metric 
of constant curvature $-1$, or equivalently, of there being a lift of the manifold 
to its universal cover, which is hyperbolic 3-space $H^3$, with the covering 
transformations acting as a discrete group of fixed point free isometries $\Gamma$.
 The manifold M is then homeomorphic to the quotient $H^3/\Gamma$.

       A hyperbolic knot $K$ in the 3-sphere $S^3$  is defined to be a knot such that
  $S^3-K$  is a hyperbolic 3-manifold. Note that the complement is a finite volume
 but noncompact hyperbolic 3-manifold.

        Given a hyperbolic knot, the SNAPPEA program, written by Jeffrey Weeks, can be
 utilized to determine the hyperbolic structure. (See ~\cite{jW03}.)In particular, the
 program yields the volume of the manifold, the symmetry group, and a variety of invariants
 that are associated to the cusps of the manifold and that will be discussed 
 in Section ~\ref{S:meridians}.
 Details of how the program determines the hyperbolic structure appear in 
 this volume (~\cite{jW04}).

     In addition, the program gives the option of determining whether or not two 
hyperbolic manifolds are isometric. Hence, since knots are known to be 
determined by their complements, (cf. ~\cite{GL89}), we can use SNAPPEA to decide 
if two given knots are the same, assuming first of all that they are both hyperbolic
 and second of all that SNAPPEA successfully finds their hyperbolic structure. In 
practice, this is one of the fastest means for determining if two hyperbolic knots
 are identical. It was utilized in the tabulation of the prime knots of 16 and 
fewer crossings(cf. ~\cite{HTW98}).

 Many knot and link complements in the 3-sphere are in fact hyperbolic. Moreover,
 knot and link complements in other 3-manifolds are often hyperbolic as well. In this 
paper, we will discuss the prevalence of hyperbolic knot and link complements and the 
various invariants associated with the hyperbolic metric carried by the complement. 
We will focus on the geometric invariants. For an excellent overview of hyperbolic 
knots with emphasis on algebraic invariants, see ~\cite{CR98}. That paper 
also contains a description of some of the applications, including the 
Smith Conjecture and the determination of symmetries of a knot.

We will not cover certain topics that are surveyed elsewhere.  See ~\cite{cG04} or ~\cite{sB02}
 for details on Dehn surgery on hyperbolic knots 
in $S^3$.  The paper ~\cite{mS98} covers unknotting tunnels and their connections with 
the hyperbolic structure on a knot complement.  Branched coverings, and 
their relationship with hyperbolicity are covered in ~\cite{CRV01} and 
~\cite{RZ01}.

 \section{What knot and link complements are known to be hyperbolic?} \label{S:what}

       The seminal work of Thurston in the 1970's and 80's demonstrated that every 
       knot in $S^3$ is either a torus knot, a satellite knot or a hyperbolic 
knot. (We include composite knots as satellite knots.)  These three categories
 are mutually exclusive.  Torus knots are well understood. Their fundamental
 groups have presentations of the form $<a,b: a^p = b^q >$. They are exactly
 the knots such that their  fundamental groups have nontrivial center. If the
 standardly embedded torus upon which a torus knot projects without crossings is
 cut open along the knot, the result is an essential annulus in the complement
 of the knot, which precludes a hyperbolic structure by fundamental group
 considerations. 

     Satellite knots have an incompressible non-boundary parallel torus in their
 complement, which also precludes a hyperbolic metric. This torus can be used to
 decompose the knot complement into simpler pieces, each of which may then be 
hyperbolic.

            The remaining case is when the knot complement has neither an 
essential annulus nor a non-boundary parallel incompressible torus. Thurston's
 remarkable theorem demonstrates that such a knot must have a hyperbolic complement. 

One might ask if it is likely that a randomly chosen knot is hyperbolic.
 If one uses a Gaussian distribution to select a cyclic sequence of n sticks glued 
end-to-end, forming a knot, then as n increases, the probability that the result
 is composite and hence non-hyperbolic, goes to one(cf.~\cite{DPS94}). If one restricts
 to prime knots, it is still the case that as n increases, the probability that the
 result is a satellite knot and hence non-hyperbolic, goes to one(cf. 
 ~\cite{dJ94}).
 So in some sense, hyperbolic knots are substantially less prevalent than non-hyperbolic
 knots. 

But on the other hand, when these non-hyperbolic knot 
 complements are cut open along essential annuli and tori, the resultant 
 pieces have a high probability of being hyperbolic. And for small 
 crossing number, the hyperbolic knots predominate. In fact, for the  2977 nontrivial  prime
 knots through twelve crossings, the only non-hyperbolic knots are seven torus knots.

     In addition, although many 3-manifolds are not hyperbolic, it was proved
 in ~\cite{rM82}  that every compact orientable 3-manifold contains a knot
 such that its complement is hyperbolic. In other words, every closed
 orientable 3-maifold is obtained by Dehn filling some hyperbolic 3-manifold.
 In this sense, hyperbolic knot complements are ubiquitous.

        Although the decomposition of the set of knots into the three classes of 
        torus, satellite and hyperbolic knots is fundamental, it does not necessarily
 allow us to easily determine
 whether a given knot is hyperbolic. It is often difficult to decide
 whether or not a given knot is a torus or satellite knot. 

       One approach is to input the knot into Jeff WeekÕs SNAPPEA program,
 which attempts to find a hyperbolic metric on the complement.  If the knot
 is hyperbolic and of reasonable crossing number, the program will almost
 always find the hyperbolic structure, thereby verifying that the knot or link
 is hyperbolic. However, if the program fails to find a hyperbolic structure,
 it could be that the knot is a torus knot or a satellite knot. Or it could be
 that the computations for determining a hyperbolic structure are too complex
 for the computer to handle with its limited memory.

      There are certain categories of knots and links in $S^3$ that are known
 to be hyperbolic. It has been proved that their complements contain no essential
 spheres, annuli or tori. We list some of these categories below:

\begin{enumerate}

\item Prime non-splittable alternating links that are not 2-braids are 
hyperbolic. 

\medskip

    This was proved in ~\cite{wM84}. This particular category of link is 
exceptionally easy to recognize as Menasco proved that an alternating link
 is splittable if and only if any 
and every alternating projection is disconnected.  He also proved that an
 alternating link is composite if and only if any and every reduced alternating
 projection has a circle that crosses the link twice and that contains crossings
 to either side. A projection is {\it reduced} if there are no crossings as in 
 Figure ~\ref{Fig:reduced}a.

\begin{figure}
\includegraphics[scale=0.5]{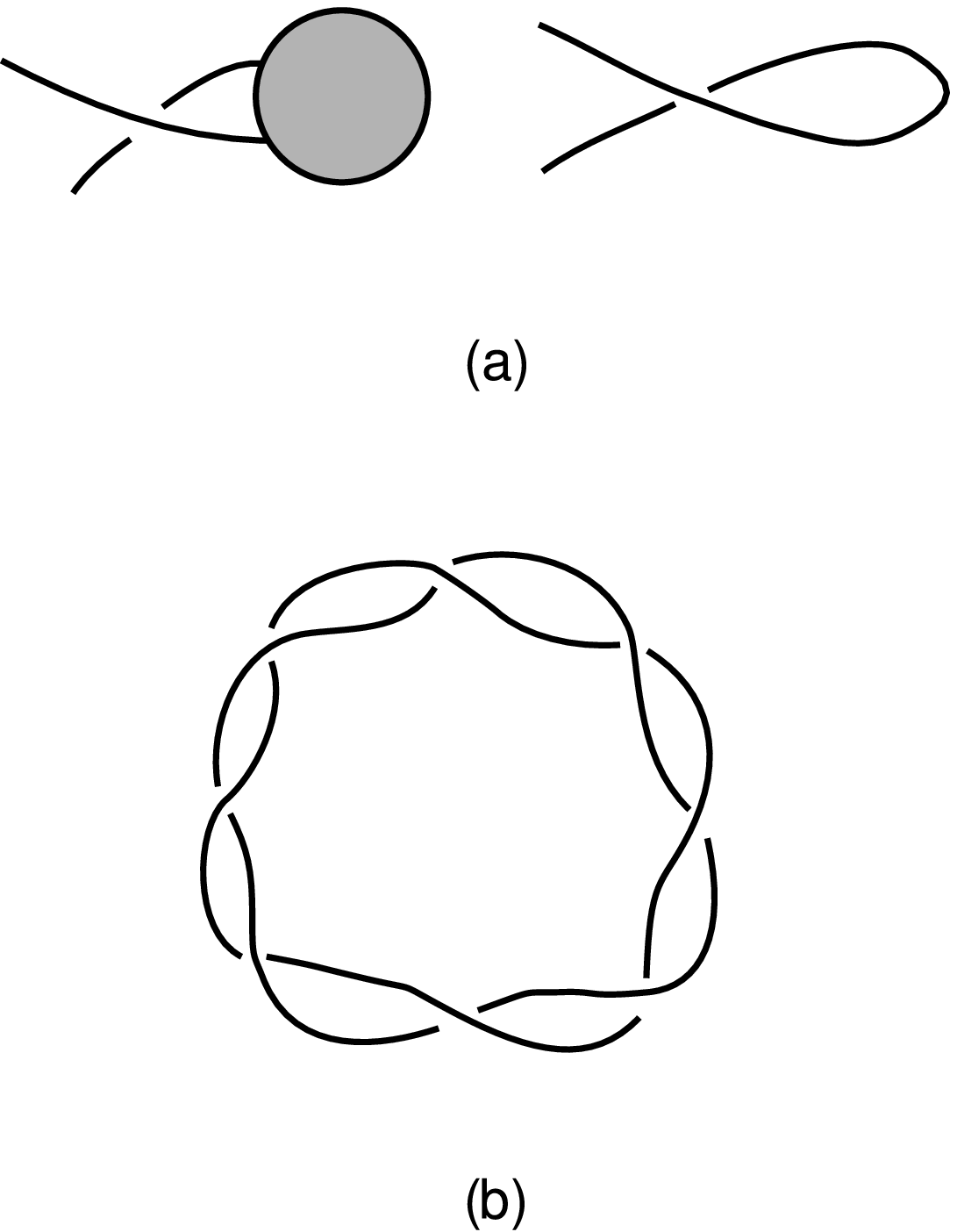}
\caption{\label{Fig:reduced} Unreduced diagrams and a 2-braid.}
\end{figure}

By a {\it 2-braid}, we mean two strands that twist around one another as in 
Figure ~\ref{Fig:reduced}b.  
If an alternating projection shows that the link is either splittable, composite 
or a 2-braid, then the complement of the link is not hyperbolic. Therefore, given
 an alternating projection of a link, we can immediately determine by examination
 whether or not it is hyperbolic.  Note that the only reduced alternating projection
 of a 2-braid is the standard one, by the Flyping Theorem of Menasco and 
Thistlethwaite(cf. ~\cite{MT91} and ~\cite{MT93}).

     Two-bridge links are all known to be prime and alternating. Hence, assuming 
that we  do not have a 2-braid (again, immediately obvious in an alternating 2-bridge
 representation), such a link is hyperbolic.

\medskip 

\item The nontrivial prime non-torus almost alternating knots are hyperbolic. (Not true for 
links.) 

\medskip 

     An almost alternating link is a non-alternating link with a projection such
 that one 
crossing change will make the projection alternating. This category of knot was
 introduced in ~\cite{Aetal92}, where it was proved that of the non-alternating 
knots of 11 or fewer crossings, all but at most three are almost alternating. 
Since then, one of those knots was shown to be almost alternating in 
~\cite{GHY01}.
  In ~\cite{Aetal92}, it was also proved that nontrivial prime non-torus almost 
alternating knots are hyperbolic. However, determining whether or not an almost
 alternating knot is prime, torus or even nontrivial remains difficult. 

\medskip 

\item Toroidally alternating knots that are prime and non-torus are 
hyperbolic.

\medskip

A toroidally alternating link is a link such that it can be projected onto the
 surface of a standardly embedded torus so that the crossings alternate over
 and under as we travel around each component on the surface of the torus and
 any nontrivial closed curve on the torus intersects the projection. An 
 example appears in 
 Figure ~\ref{Fig:toralt}. Almost 
alternating knots are all toroidally alternating. In ~\cite{cA94}, it was proved that
 nontrivial prime non-torus toroidally alternating knots are hyperbolic. But 
again, showing that a toroidally alternating knot is not trivial, composite
 or torus remains a difficult question. 

\begin{figure}
\includegraphics[scale=0.7]{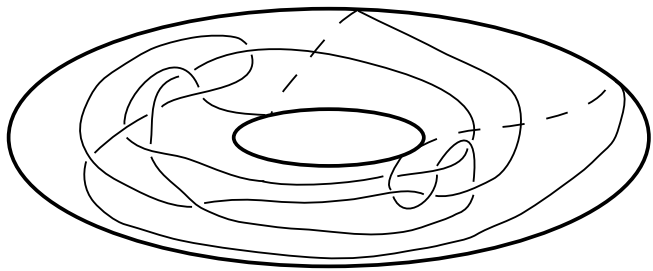}
\caption{\label{Fig:toralt} A toroidally alternating knot.}
\end{figure}

\medskip 

\item Augmented alternating links are almost all hyperbolic. 

\medskip

Given a prime alternating link projection, one can add trivial components that 
bound disjoint disks, each of which is perpendicular to the projection plane 
and  intersects the knot at exactly two points. If the resulting link is 
nonsplittable, and no two of these added components are parallel in the 
resulting link complement, we call the result an augmented alternating 
link. See  Figure ~\ref{Fig:augalt}.

\begin{figure}
\includegraphics[scale=0.5]{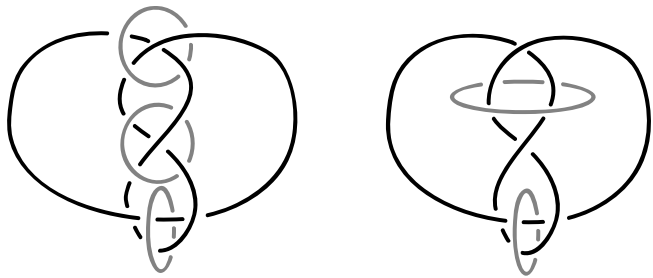}
\caption{\label{Fig:augalt} Two augmented alternating links from the 
figure-eight knot.}
\end{figure}

     In ~\cite{cA86}, it was proved that if the initial link is not a 2-braid,
 then an augmented alternating link is hyperbolic. Note that in 
 ~\cite{mL03},
 Lackenby proved that the closure in the geometric topology of the set of all
 hyperbolic alternating links
 is the set of all hyperbolic alternating links and augmented alternating links.

      In addition, one can cut open along a twice-punctured disk, twist any number
 of half-twists and reglue. If the number of half-twists is even, we 
 simply obtain a new link whose complement is homeomorphic to the original.
 If the number of 
 half-twists is odd, we obtain a new link complement. In ~\cite{cA85}, 
it was proved that if the original link
 complement is hyperbolic, so is the new link complement, and it will have
 the same volume.

    So, given a minimal projection of any nontrivial nonsplit link complement other
 than the standard projection of the two-braid or an obviously composite projection,
 one can add these trivial components around crossings, except where doing so yields
   parallel copies, to obtain a link, which is hyperbolic, since we can switch the
 crossings on the original link to make it alternating. We call such a link fully 
augmented. See ~\cite{eS03} for more details.

\medskip 

\item Arithmetic Links

\medskip

Some of the first link complements known to be hyperbolic were arithmetic, meaning
 that their fundamental groups are commensurable 
 (up to conjugacy in ${\rm PSL}(2,\bold C)$) with a Bianchi group ${\rm PSL}(2,O_d)$,
 where $d$ is a square-free positive integer, and  $O_d$ is the ring of integers in
 the imaginary quadratic field $\Q(\sqrt{-d})$.
 See for instance, ~\cite{rR75,rR79,wT78, nW78, aH83, mB92, jS96, jS99, mB01}.
 Although it appears that arithmetic links are a small
 subset of the set of all links, Baker proved in ~\cite{mB02} that every link in $S^3$ is a 
sub-link of an arithmetic link in $S^3$.
 In the case of knots, Reid proved that the only arithmetic knot is the 
 figure-eight knot(cf. ~\cite{aR191}).

\medskip 

\item Montesinos links 

\medskip

A Montesinos link (also called a star link) is obtained by conecting n 
rational tangles in a simple cyclic fashion. Such a knot or link is denoted
 $K(\frac{p_1}{q_1},\dots, \frac{p_n}{q_n})$, where $\frac{p_i}{q_i}$ 
 denotes the ith rational tangle. 
 In ~\cite{uO84}, it was  proved that a
 Montesinos link  $K(\frac{p_1}{q_1},\dots, \frac{p_n}{q_n})$ with $q_i \ge 2$, is 
hyperbolic if it is not a torus link and  not equivalent to 
$K(1/2,1/2, -1/2, -1/2)$, $K(2/3, -1/3, -1/3)$,
 $K(1/2, -1/4, -1/4)$, $K(1/2,-1/3, -1/6)$ or the mirror image of these links. In 
~\cite{BS80} the torus links that are Montesinos links  are identified.

\medskip 

\item Mutants of a hyperbolic link.

\medskip

Given a knot or link in a projection and a circle in the projection plane that intersects
 the knot at four points and separates the knot into two tangles,
 one can perform the following operation. Cut the knot open at these four points and flip
 the interior tangle either around a vertical or horizontal axis, or rotate it 180 
degrees about an axis perpendicular to the projection plane before reattaching it
 to the outside tangle at the four points. This operation is called 
 {\it mutation} and the resulting knots are called {\it mutants} of the 
 original. In ~\cite{dR87}, Ruberman demonstrated that a mutant of a hyperbolic
 knot or link is also hyperbolic and it has the same volume.

\medskip 

\item Belted sums of hyperbolic links.

\medskip

Given two links, each with a trivial component bounding a disk punctured twice by the 
link, one can cut them open along the disks in each, glue the resultant 
disks from the one complement to the disks from the other appropriately, to 
obtain a link complement as in Figure ~\ref{Fig:belt}. In ~\cite{cA85}, it was proved
 that the belted sum of hyperbolic links is always
 hyperbolic with volume equal to the sum of the volumes of the original links.

\begin{figure}
\includegraphics[scale=0.5]{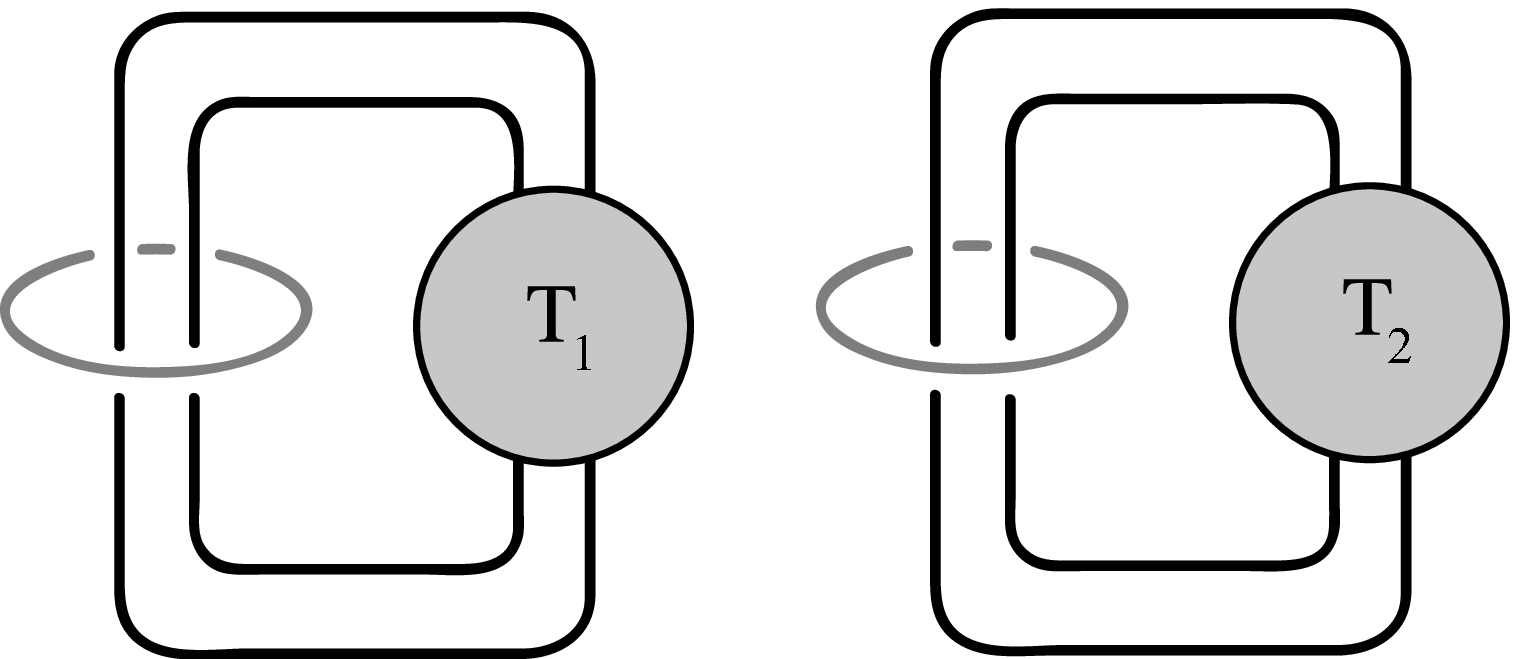}
\includegraphics[scale=0.5]{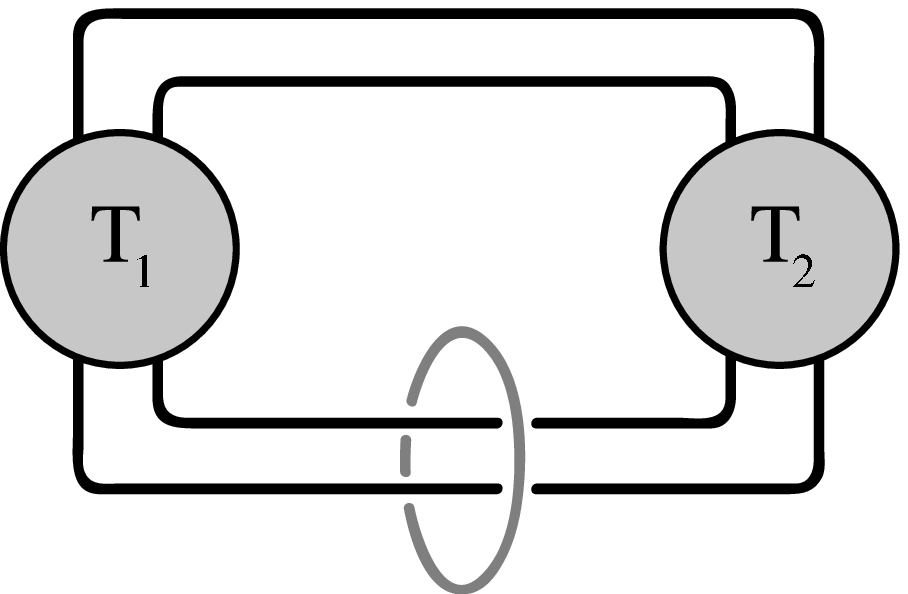}
\caption{\label{Fig:belt} The belted sum of two links.}
\end{figure}

\medskip 

\item Hyperbolic Tangles.

\medskip

 In  ~\cite{yW96}, Wu determined all of the non-hyperbolic 
alegbraic tangles, finding them to be a small subset of the set of all 
algebraic tangles.  In particular, the hyperbolic algebraic tangles can be closed off to
 obtain algebraic 
knots and links that are hyperbolic. This gives a broad set of examples. 
Note that the tangles themselves can be thought of as complements of genus 
two handlebodies, and as such, are realized as hyperbolic sub-manifolds of 
$S^3$ with totally geodesic boundary.

In ~\cite{kS99}, it is shown that certain n-string alternating tangles 
are hyperbolic. These tangles can be utilized to show that prime 
semi-alternating links are hyperbolic. (See ~\cite {LT88} for the definition 
of semi-alternating links.) 

\end{enumerate}

\section{Volumes of Knots}\label{S:volumes}

The most natural invariant associated to hyperbolic knots is the volume
 of the complement.  Work of Thurston and J\o rgensen shows that knot volumes
 are a well-ordered subset of $\R$. In ~\cite{AHW91}, lists of volumes of 
 knots through ten crossings were calculated. Cao and Meyerhoff (cf. 
 ~\cite{CM01}) proved
 that the figure-eight knot has the smallest possible volume for a knot 
 complement, the volume of which is $2.02988\dots= 2v_3$, where $v_3 = 
 1.01494\dots$ is the volume of an ideal regular tetrahedron. In fact, they 
 proved
 that this is the smallest volume of any noncompact orientable hyperbolic 
3-manifold, the volume being shared by the figure-eight knot complement and
 one other manifold known as the sibling of the figure-eight knot complement.

A hyperbolic n-component link is known to have volume at least $nv_3$ 
(cf.~\cite{cA88}), and a 
2-component link is known to have volume at least $2.3952 
v_3$(cf.~\cite{hY01}), although 
the expectation is that this bound is not sharp. The smallest known 
2-component link is the Whitehead link, with volume $3.6638\dots$.

    Results of Thurston and J\o rgensen demonstrate that if one does
  $(p,q)$-Dehn filling on a hyperbolic knot or link complement, with
 $p^2 + q^2$ large enough, the resulting manifold will also be hyperbolic
 with volume less than the volume of the original manifold but approaching
  the volume of the original manifold as $p^2 + q^2$ approaches 
  $\infty$. (See ~\cite{wT78}.)

\begin{figure}
\includegraphics[scale=0.5]{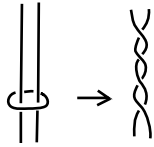}
\caption{\label{Fig:twist} Twisting by surgery.}
\end{figure}

      In particular, we can take a component as in Figure ~\ref{Fig:twist}, and do
 $(1,p)$-surgery to obtain a knot with arbitrarily many crossings but
 volume bounded by the original manifold. 

\begin{figure}
\includegraphics[scale=0.5]{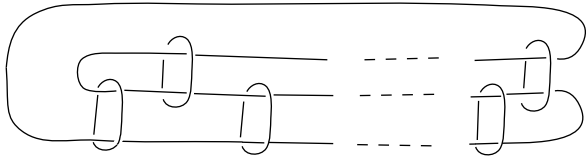}
\caption{\label{Fig:twobridge} Two-bridge knots come from surgery on these 
links.}
\end{figure}

     On the other hand, there are knot complements with arbitrarily
 large volume. We will show this is true for the two-bridge knots. 
Two-bridge knots are obtained from surgery on the smaller 
components in the links depicted in Figure ~\ref{Fig:twobridge}.
 Each of these links is a belted sum (see the previous section for the
 definition and properties) of Borromean rings. 

      The belted sum of n copies of the Borromean rings has volume equal
 to $n (7.3276\dots)$. When we do high surgery on the augmenting components,
 we obtain two-bridge
 knots with volume arbitrarily close to that of the belted sum. Hence, volumes
  of knot complements can be arbitrarily high.

      How effective an invariant is volume for distinguishing between knots?
 In general, it is very good. The first example of two knots with the same
 volume doesn't occur until we consider knots of up to twelve crossings. The
  $5_2$ knot and a 12-crossing knot do have the same volume. 

        However, the operation of mutation preserves volume. Hence, many
 knots do share their volumes with other knots. Moreover, the operation
 of twisting along a twice-punctured disk preserves volume as well,
 yielding a variety of link complements
 with the same volume.

     Can we say anything about the volume of a particular knot or link 
complement simply by looking at its projection? In ~\cite{cA83}, it was 
proved that the volume of the complement of an n-crossing hyperbolic knot other than  
the figure-eight knot is bounded above by $(4n-16)v_3$.

 In ~\cite{mL03}, 
Lackenby defines a twist in an alternating projection to be a maximal chain of
 adjacent bigon 
regions(as in the second part of Figure ~\ref{Fig:twist}),
 or to be a single crossing that is not adjacent to a bigon region. The twist
 number of a projection is the number of twists within it. He proves(with an
 improvement of his upper bound due to Ian Agol and 
Dylan Thurston) that a hyperbolic alternating knot or link in a prime 
alternating projection of twist number $t$ satisfies 
$$v_3 (\frac{t-2}{2}) \le Volume(S^3-K) < v_3 (10t-16)$$  By choosing a twist 
 reduced projection, which always exists, the lower bound can be 
improved to $v_3(t-2)$.

Is hyperbolic volume related to the more recently defined quantum 
invariants?  In ~\cite{rK97}, Kashaev conjectured that the hyperbolic volume of a
 knot complement is determined by the asymptotic behavior of a link invariant
 that depends on the quantum dilogarithm and that was introduced by Kashaev 
 in ~\cite{rK95}. This is known as the Kashaev Conjecture. It
 has been verified for a handful of knots.

\section{Cusps}\label{S:cusps}

\begin{figure}
\includegraphics[scale=0.5]{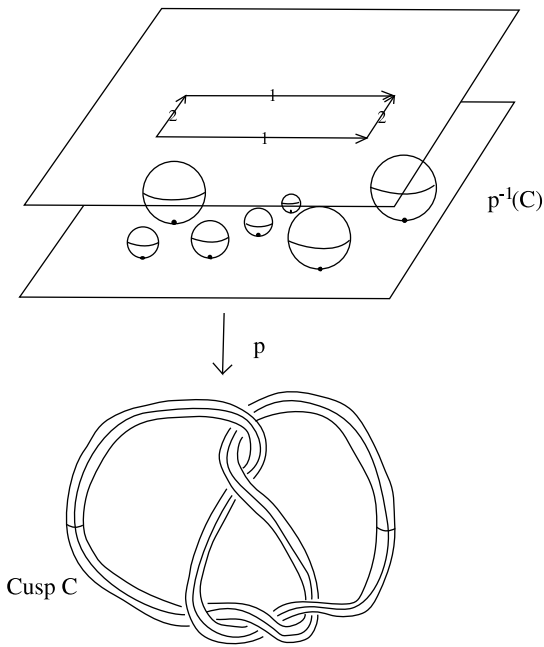}
\caption{\label{Fig:cusp} The cusp of a knot.}
\end{figure}

Given a hyperbolic knot in $S^3$, one can define a cusp $C$ for the knot
 to be a neighborhood of the missing knot such that it lifts to a set
 of horoballs with disjoint interiors in the universal cover $H^3$. See 
 Figure ~\ref{Fig:cusp}.

Topologically, a cusp is homeomorphic to $T^2 \times [0,1)$ where 
$T^2 \times 0$ corresponds to the boundary of the neighborhood and
 the missing $T^2 \times 1$ corresponding to the knot itself. Choosing
 one of the covering horoballs to be the horoball centered at {$\infty$}
 in the upper-half-space model of $H^3$, the subset of covering translations
 sending this horoball back to itself are all Euclidean translations,
 generated by two translations. A fundamental domain for their action
 is a parallelogram P in the horizontal boundary plane and the collection of horizontal
 parallelograms directly above this one.  The parallelogram P projects
 to the torus boundary of the cusp. As we move up the z-axis in hyperbolic
 space, each horizontal parallelogram projects to a concentric torus
 around the missing knot, the set of them shrinking in toward the missing
 knot. 

     In $H^3 \cup \partial H^3$, the Euclidean translations fix a single
 point on the boundary, that being the point at $\infty$. We call an isometry
 that fixes a single point on the boundary a {\it parabolic} isometry.  All such
 isometries will send the horoballs tangent to the boundary at that fixed
 point back to themselves. All other isometries in the group $\Gamma$ of covering
 translations have two fixed points on the boundary. These are called
 {\it hyperbolic} isometries. They correspond to translations along the geodesic
 with the two fixed points as endpoints, together with a possible rotation
 about the geodesic.

    The third type of orientation preserving isometry is a pure rotation about
 a geodesic. But as it is not fixed point free in $H^3$, it cannot appear as
 a covering translation in $\Gamma$.
     In the case of a knot, we have a single cusp in the complement. We can expand
 that cusp until it touches itself. Then, the set of horoballs that are the pre-image
 of $C$ in $H^3$ will expand until two first touch. Since every horoball in 
 $p^{-1}(C)$
 is identified to every other one, the horizontal plane that is the boundary
 of the horoball centered at $\infty$ will touch other horoballs. We call this
 cusp a {\it maximal cusp}. In the case we have a link, we can expand the cusps until
 they touch each other or themselves, to obtain a maximal set of cusps. However,
 in different situations, it may be appropriate to do the expansions in different
 ways. For instance, we may insist that all the cusps have the same volume
 or the same length of shortest curve in their boundaries while we expand the
 set. Or we may choose to put the largest possible volume in the expanded set
 of cusps.

      We define the {\it cusp density} of a hyperbolic knot or link complement to be
 the ratio of the largest possible volume in a maximal set of cusps divided by
 the total volume. Meyerhoff noted in ~\cite{rM86} that this number is at most $0.853É =
 \frac{\sqrt{3}}{2 v_3}$ where again, $v_3 = 1.01494\dots$ is the volume of an ideal regular
 tetrahedron. The cusp density of the figure-eight knot complement realizes
 this upper bound. In ~\cite{Aetal02}, it was proved that there are knots of arbitrarily
 small cusp density, although other methods for generating such knots were
 previously known, but had not appeared in the literature. Although it is known
 that the set of cusp densities for hyperbolic 3-manifolds are dense in the
 interval $[0, .853\dots]$, it is not known whether the same holds true for cusp
 densities of  hyperbolic link complements or perhaps, even  for cusp
 densities of hyperbolic knot complements.

\section{Meridians and other cusp invariants}\label{S:meridians}

The Thurston-Gromov $2\pi$ theorem implies that for any Dehn filling along
 a curve $c$ in the cusp boundary of a hyperbolic 3-manifold, the resultant
 manifold will be negatively curved if  $c$ has length at least $2\pi$. It 
 is conjectured that in fact the resultant manifolds are hyperbolic with 
 constant negative curvature, but this question remains open.
(See ~\cite{BH96} for more details.) If we define the length $|m|$ of a
 meridian for a hyperbolic knot K to be the length of its shortest representative
 in the maximal cusp boundary, then $|m|< 2\pi$, since Dehn filling along
 the meridian yields $S^3$, which is not negatively curved.
 
     In ~\cite{mL00} and ~\cite{iA00}, a variation on the $2\pi$ Theorem, 
     known as the 6-Theorem, was proved. It shows that a manifold obtained by Dehn filling along
 a minimal length simple closed curve of length greater than 6 in the cusp boundary of a
 hyperbolic manifold is {\it hyperbolike}, which is to say that it is irreducible, boundary-irreducible, 
 and has infinite fundamental group that is word hyperbolic. In other 
 words, it has these attributes that one would expect of a negatively curved 
 manifold. Agol also gave 
 an example of a cusped manifold such that there was a curve of length 6 
 in the cusp boundary such that surgery on that curve yielded a 
 non-hyperbolike manifold. So the 6 bound is sharp. In ~\cite{Aetal03}, we produce 
 examples of knot complements with longitude of length exactly 6 such that 
 surgery on the longitude yields non-hyperbolike manifolds. So the 
 6-theorem is sharp for knot complements in $S^3$.

Since Dehn filling
 along the meridian of a knot yields the 3-sphere, which has trivial fundamental
 group, the 6-Theorem implies every meridian of a knot must have length no greater than 6. 
Moreover, any nontrivial curve in the cusp boundary must have length at least 1.
 This is because a maximal cusp when lifted to the horoball centered at 
$\infty$, will have other horoballs tangent to it. The shortest parabolic
 isometry must shift the horoballs at least a distance 1 so that they do 
 not overlap with one another. 

     So all meridian lengths fall in the interval [1, 6]. In ~\cite{cA202},
 it was proved that there is only one knot that realizes the lower bound.
 The length of the meridian in the figure-eight knot complement is exactly 1,
 and no other hyperbolic manifold has a nontrivial curve in its maximal
 cusp boundary this small. Moreover, in ~\cite{cA03}, it is proved that
 the next shortest meridian is that of the $5_2$ knot, with a length of
 $1.150964\dots$, and that there are no other knots of meridian length
 less than $\sqrt[4]{2}= 1.189207\dots$.

\begin{figure}
\includegraphics[scale=0.5]{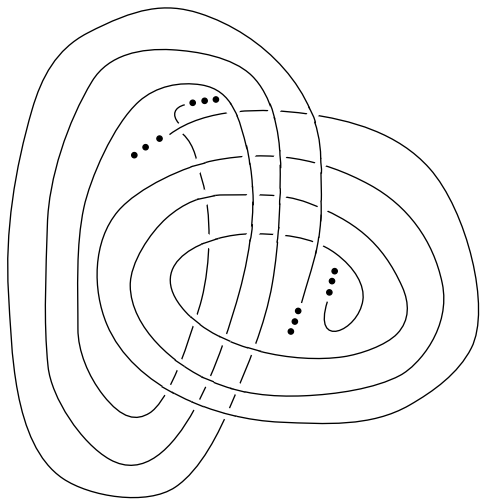}
\caption{\label{Fig:agolknot}These knots have meridian length approaching 
4 from below.}
\end{figure}

     What about the upper bound? In ~\cite{iA00}, Agol gives an example
 of a sequence of knots with meridian lengths approaching 4 from below. We 
 show how to construct such a sequence of knots in Figure ~\ref{Fig:agolknot}. 
 As the number of times that the knot wraps around itself both vertically 
 and horizontally increases, the meridian length approaches 4 from below. These
 are the largest known meridian lengths to date. 

   Utilizing ideas from ~\cite{zH98}, it was proved in ~\cite{Aetal02} that the meridian
 length of a knot $K$ is always at most $6-\frac{7}{c}$, where $c$ is the crossing 
 number of $K$.
 The idea is to take the singular punctured surface obtained by coning the
 knot to a point below the projection plane. If it is not incompressible
 or boundary-incompressible, then we can do compressions and boundary-compressions
 that will only improve the ultimate bounds. Then we pleat the surface
 and it inherits a 2-dimensional hyperbolic metric from the 3-dimensional
 hyperbolic metric on the manifold. The length of the boundary curves of
 the surface intersected with the cusp is equal to the area of the intersection
 of the  cusp with the surface. This is bounded above by $\frac{3}{\pi}$ times the area
 of the surface, which is given by $2\pi\abs{\chi(S)}$. From these inequalities,
 we can obtain bounds on meridian length and cusp volume.

     One can obtain substantially better bounds for certain classes of knots.
 In ~\cite{cA96}, it was proved that the two-bridge knots have meridian lengths in the
 interval [1,2), with sequences of meridians approaching 2.    

In ~\cite{Aetal02}, upper bounds for alternating knots were determined. In particular,
 it was shown that the meridian of an alternating knot with crossing number
 $c$ is bounded above by $3-\frac{6}{c}$. The proof utilizes the checkerboard surfaces
 that come from any given projection of a knot. In the case of alternating
 knots, it was proved in ~\cite{MT93} that these surfaces are incompressible and
 boundary-incompressible. Pleating these surfaces allows one to to prove 
 the result.

The expectation is that the actual upper bound for the meridian length 
of alternating knots is substantially smaller.
 For example, we say a property J is true for almost all hyperbolic 3-manifolds
 if for every fixed volume $V > 0$, J is true for all but finitely many 
hyperbolic 3-manifolds with volume less than V.  In ~\cite{eS03},it is  
proved that almost all alternating hyperbolic knots have
 meridian length bounded above by 2. The expectation is that 2 is the correct 
 bound for all alternating hyperbolic knots.
    
In addition, in ~\cite{eS03}, Schoenfeld proves that if we
 fully augment any reduced projection of a nontrivial knot other than a 
 2-braid, as in Section 
~\ref{S:what}, the meridian of the original knot becomes 2 in the resultant 
link complement. In particular, 
this means that if we take any reduced projection of any non-2-braid knot
 and twist about each of its crossings to add bigons, the resultant knots
 will have meridian approaching 2. 

In ~\cite{zH98}, an application of meridian length to the determination of 
crossing number of a knot is given. Specifically, it is shown that if $C$ is a maximal
cusp for a hyperbolic knot $K$, then the crossing number ${\rm
cr}(K)$ satisfies $${\rm cr}(K)\geq{\rm area}(\partial C)/(|m|(2\pi-|m|)).$$ 
Moreover, if $K'$ is a satellite of the hyperbolic knot $K$ of degree
$p$, then ${\rm cr}(K')\geq p^2{\rm area}(\partial C)/(|m|(2\pi-|m|))$.

 \section{Geodesics and Totally Geodesic Surfaces}\label{S:geodesics}

Geodesics in hyperbolic knot complements may or may not be closed. If  a 
geodesic is closed, it might intersect itself or not. However, every hyperbolic knot 
or link complement contains a simple closed geodesic(~\cite{AHS99}). 

A {\it systole} is a shortest closed geodesic, and the {\it systole length} is the 
length of the shortest closed geodesic in a manifold. If we take two 
strands in a knot complement, and obtain a sequence of knots by twisting 
these two strands around one another as in  Figure ~\ref{Fig:twist}, 
then the length of the 
geodesic that wraps around the two strands will be approaching 0. Hence, 
this yields a sequence of knot complements with systole length shrinking to 
0. So there is no lower bound for systole length in knot and link 
complements.
        On the other hand, there is an upper bound. Although the systole 
length can be arbitrarily large for hyperbolic 3-manifolds in general,
 it was proved 
in ~\cite{AR02} that systole length for hyperbolic knots and links in $S^3$ 
is bounded above by 7.35534. For hyperbolic alternating knots, this was 
improved to 4.5 in ~\cite{Aetal02}. It would be interesting to obtain 
better bounds on systole length for a variety of categories of knots and links.

Closed geodesics in knot complements can themselves be either unknotted or  knotted as curves 
in $S^3$. See for instance, ~\cite{sM01}, where  a 
variety of knotted geodesics in the figure-eight knot complement are 
displayed.

In ~\cite{tS91} (see also~\cite{sK88}), it was proved that the complement
 of a simple closed geodesic or a set of 
disjoint simple closed geodesics in a hyperbolic manifold will itself 
yield a hyperbolic manifold. One can look at the link complements obtained from 
removing certain geodesics from a knot or link complement, in the hopes of 
decomposing the resultant manifolds into certain canonical pieces.

One might hope that the complement of a simple closed geodesic in a 
hyperbolic knot complement has minimal volume among all the complements 
of simple closed curves in the same free homotopy class. But counterexamples to 
this conjecture are given in the figure-eight knot complement in ~\cite{sM01}.

\medskip

We turn now to surfaces in knot complements. A closed embedded incompressible
 surface in a hyperbolic knot complement can
 come in one of two varieties.  The first possibility is that it is 
quasi-Fuchsian. This means that it lifts to the disjoint union of 
topological planes in $H^3$, each with a limit set that is a quasi-circle.
 In the case that the limit set is an actual circle, and  the planes are geodesic,
 and we say that the surface is {\it Fuchsian} or {\it totally geodesic}. 

 The second possibility for a closed embedded incompressible surface is
 that there are simple closed curves on the surface that can be homotoped
 through the knot complement into a neighborhood of the missing knot.
 This means that the corresponding isometry, when we lift to hyperbolic
 space, is a parabolic isometry. We call such a curve an {\it accidental
 parabolic curve} and we call such a surface an {\it accidental surface}.  In 
 ~\cite{IO200},
 it was proved that if an accidental parabolic curve exists for a
 given surface $S$, then it is isotopic to a unique curve on the boundary
 torus of the cusp. In particular, there is a well-defined
 accidental slope for each accidental surface. However, there can be
 more than one accidental slope for a given knot. It follows from
 ~\cite{CGLS87} that accidental slopes must be meridional or integer.
 
    In ~\cite{AR93}, explicit examples were given of closed quasi-Fuchsian surfaces
 in knot complements. These surfaces were shown not to be totally geodesic.
 In fact, in ~\cite{MR92}, it was conjectured that hyperbolic knot complements
 in $S^3$ do not contain any closed embedded totally geodesic surfaces. 
This has been proved for hyperbolic knots that are alternating 
knots, tunnel number one knots, 2-generator knots and knots of braid index
 three~\cite{MR92}, almost alternating knots ~\cite{Aetal92},  toroidally
 alternating knots ~\cite{cA94}, Montesinos knots ~\cite{uO84}, 3-bridge 
 knots and double torus
 knots ~\cite{IO200}) and knots of braid index three  ~\cite{LP85} and  four ~\cite{hM02}.
Note that the conjecture does not hold for links. An explicit
 counterexample is given in ~\cite{MR92}.

There are closed totally geodesic surfaces immersed in hyperbolic knot 
complements. In ~\cite{aR291}, Reid shows that the figure-eight knot 
complement contains infinitely many such non-homotopic surfaces. See 
~\cite{AR97} for two additional knots with immersed totally geodesic 
surfaces in their complements.

   An incompressible boundary incompressible surface $S$ with boundary
 properly embedded in a hyperbolic knot exterior can have one of three possible behaviors:

\begin{enumerate}

\item $S$ can be quasi-Fuchsian.
\item $S$ can be accidental.
\item $S$ can be a virtual fiber in a fibered knot, with limit set the entire
 boundary of $H^3$.

\end{enumerate}
     
     Specific examples of incompressible boundary-incompressible surfaces
 are afforded by minimal genus Seifert surfaces.  In ~\cite{sF98}, it is proved
 that a minimal genus Seifert surface in a non-fibered hyperbolic knot
 complement must be quasi-Fuchsian. It can never be accidental. 

\begin{figure}
\includegraphics[scale=0.5]{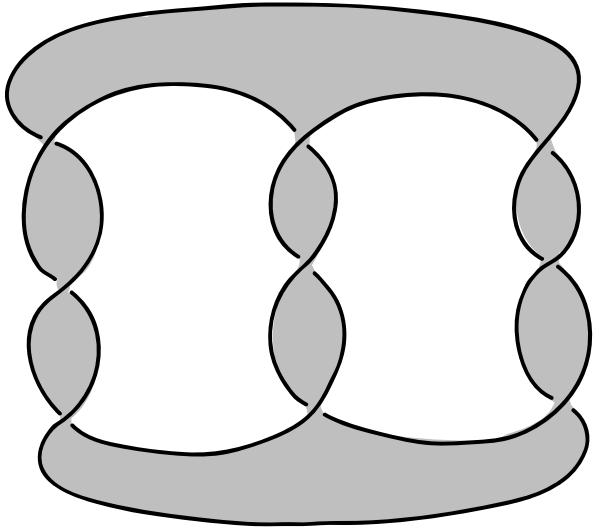}
\caption{\label{Fig:pretzel}A (3,3,3)-pretzel knot has a totally geodesic Seifert surface.}
\end{figure}

     Knots can have totally geodesic minimal genus Seifert surfaces. In
 ~\cite{AS03}, examples such as a (p,p,p) pretzel knot(Montesinos knot 
 $K(1/p,1/p,1/p)$) are shown to have totally
 geodesic Seifert surfaces.(See Figure ~\ref{Fig:pretzel}.) There are examples of hyperbolic knots with
 totally geodesic Seifert surfaces, where the Seifert surfaces are both
 free (the complement of a neighborhood of the Seifert surface in the knot 
 complement is a handlebody) and non-free.

      But the expectation is that  knots with totally geodesic Seifert
 surfaces are the exception. For instance, in  ~\cite{AS03}, it is proved that
 2-bridge knots never have totally geodesic Seifert surfaces.

An interesting question is whether there are totally geodesic
 surfaces in knot complements other than such Seifert surfaces. This is 
 one of many open questions that still remain in the theory of hyperbolic 
 knots.


\end{document}